\documentclass[a4paper, 12pt]{article}
\usepackage{geometry}
\usepackage[english]{babel}
\usepackage[utf8]{inputenc}
\usepackage{amsmath}
\usepackage{commath}
\usepackage{amsthm}
\usepackage{amssymb}
\usepackage{setspace}
\usepackage{mathtools}
\usepackage{color}

\usepackage[
style=alphabetic,
sorting=ynt,
doi=false,isbn=false,url=false,eprint=false
]{biblatex}
\def\bfa{{\bf a}}

\def\bfh{{\bf h}}

\def\bfv{{\bf v}}
\def\bfw{{\bf w}}

\def\bfx{{\bf x}}
\def\bfy{{\bf y}}
\def\bfz{{\bf z}}
\def\det{{\rm det}}

\newtheorem{theorem}{Theorem}[section]

\newtheorem{lemma}[theorem]{Lemma}
\newtheorem{example}[theorem]{Example}
\newtheorem{definition}[theorem]{Definition}
\newtheorem*{remark}{Remark}

\title{The Harmonic GBC Function Map is a Bijection if the Target Domain is Convex }
\author{Chongyang Deng\footnote{School of Science, Hangzhou Dianzi  
University, Xiasha, Hangzhou 310018, China. This author is supported 
by  the National Science Foundation of
China under grant \# 61872121 and Zhejiang Provincial Science and 
Technology Program in China under Grant \# 2021C01108.}  
\and Tsung-wei Hu\footnote{Department of Mathematics,
University of Georgia, Athens, GA 30602.} 
\and Ming-Jun Lai 
\footnote{mjlai@uga.edu. Department of Mathematics,
University of  Georgia, Athens, GA 30602. 
This author is supported by the Simons Foundation collaboration grant \#864439.} }
\begin{document}
\maketitle

\begin{abstract}
Harmonic generalized barycentric coordinates (GBC) functions have been used for cartoon 
animation since an early work in 2006\cite{JMDGS06}. 
A computational procedure was further developed in \cite{SH15} for 
deformation between any two polygons. 
The bijectivity of the map based on harmonic GBC functions is still murky in the literature.   
In this paper we present an elementary proof of the bijection of the  harmonic 
GBC map transforming from one arbitrary polygonal domain $V$ to a convex polygonal domain $W$.  
This result is further extended to a more general harmonic map from one simply connected 
domain $V$ to a convex domain $W$ 
if the harmonic map preserves the orientation of the boundary of the  domain $V$. 
In addition, we shall point out that the harmonic GBC map is also a diffeomorphism over the 
interior of $V$ to the interior of $W$.    
Finally,  we remark on how to construct a harmonic GBC map from $V$ to $W$ 
when  the number of vertices of $V$ is different from the number of vertices of $W$ and how to
construct harmonic GBC functions over a polygonal domain with a hole or holes.     
We also point out that it is possible to use the harmonic GBC map to 
deform a nonconvex polygon $V$ to another nonconvex polygon $W$ by a 
good arrangement of the boundary map between $\partial V$ and $\partial W$. 
Several numerical deformations based on images 
are presented to show the effectiveness of the map based on bivariate spline approximation 
of the harmonic GBC functions.    
\end{abstract}

\bigskip

\leftline{{\bf Mathematics Subject Classification:} 31A05, 35J25, 30C60, 53A10}

\section{Introduction}
We are interested in the bijective property of the map based on harmonic GBC functions 
which transforms from one arbitrary 
polygonal domain  $V \subset \mathbb{R}^2$ to a convex polygonal domain 
$W\subset \mathbb{R}^2$. Let us first recall the GBC functions.  
Given a polygon $V \in \mathbb{R}^2$ of $n$-sides with vertices $\bfv_i, i=1, \cdots, n$,
the functions $\phi_i, i=1, \cdots, n$ satisfying the following conditions:  
\begin{equation}
\label{GBC}
\begin{cases}
\sum_{i=1}^n \phi_i(\bfx)  &=  1 \cr
\sum_{i=1}^n \phi_i(\bfx)\bfv_i &=  \bfx\cr
\phi_i(\bfx) &\ge 0, \quad i=1, \cdots, n,
\end{cases}
\end{equation}
for all $\bfx=(x,y)\in V$ and being linear on edges of $V$ 
with $\phi_i(\bfv_j)=1$ if $i=j$ and $=0$ if $i\not=j$. 
Such functions are called generalized barycentric coordinates (GBC) or GBC functions over $V$.  
 
The study of generalized barycentric coordinates (GBC) started from a seminal work 
in \cite{W75}. Since then, there are many GBCs which have been constructed. For example, 
harmonic GBC functions in the 
3D setting was explained in \cite{JMDGS06} together with many interesting and useful properties. 
We refer to  a recent survey in \cite{F15} and a book \cite{HS18} edited by leading experts
K. Hormann and N. Sukumar. Recently, a general minimization method is proposed to construct 
GBC functions (cf. \cite{DFL20}). 
More precisely, consider the following minimization problem:
\begin{equation}
\label{min}
\begin{array}{ll}
\min\{  F(\phi_1, \cdots, \phi_n): \hbox{subject to }& \sum_{i=1}^n \phi_i(\bfx)=1, 
 \sum_{i=1}^n\phi_i(\bfx) \bfv_i =\bfx, \cr 
& \bfx\in V, \phi_i\ge 0, i=1, \cdots, n\}, 
\end{array}
\end{equation}
where $F$ is a convex function and $V$ has $n$ boundary vertices. For example, one can choose 
$$\displaystyle F(\phi_1, \cdots, \phi_n)=\sum_{i=1}^n \int_{V}|\nabla \phi_i(\bfx)|^2 d\bfx$$ 
which leads to a new way to compute harmonic GBCs based on bivariate spline functions 
(cf. \cite{ALW06} and \cite{LS07}). 
We refer to \cite{DFL20} for a detail explanation.     
The GBCs have found their applications in 
geometric design. See, e.g. \cite{JBPS11}, \cite{Zetc14}, \cite{LL21}, and etc.. 
In addition, they found their applications in 
numerical solution of partial differential equations 
(e.g. \cite{MRS14}, \cite{RGB14}, \cite{FL16},  and etc.). 

One interesting property of the GBCs functions is that 
it can form a transform or map from one polygon to another polygon
easily. That is, letting $V$ be a polygon with boundary vertices $\bfv_i, i=1, \cdots, n$ 
arranged in the counter-clockwise direction and $W$ be another polygon with boundary vertices 
$\bfw_i, i=1, \cdots, n$ in the counter-clockwise direction, 
if $\phi_i, i=1, \cdots, n$ are GBC functions over $V$ satisfying  (\ref{GBC}), we can define
\begin{equation}
\label{GBCmap}
{\bf F}(x,y)= \sum_{i=1}^n \phi_i(x,y)\bfw_i, \quad (x,y)\in V.
\end{equation}
Then ${\bf F}$ is a map from $V$ to $W$ when $W$ is a convex polygon.  Clearly, the image 
$Im({\bf F})\subset W$ due to the properties of (\ref{GBC}). One can show that the map 
is onto which will be proved in the next section.   
The injectivity of the Wachspress  from a convex polygon $V$ to a 
convex polygon $W$ was established in \cite{FK10}. 
However, the mean value coordinates may not be injective when 
a convex polygon $V$ is a pentagon as shown 
in \cite{FK10}. It is interesting to know if the harmonic GBC  map ${\bf F}$ 
is injective or not. 
This paper is also motivated by the work in \cite{SH15}, 
where the researchers proposed an approach to construct a map ${\bf F}$ from a polygonal 
domain $V$ to another polygonal domain $W$ by introducing an intermediate 
convex polygonal domain $\Theta$ such that 
$$
{\bf F}= \Phi_1^{-1}\circ \Phi_0, 
$$ 
where $\Phi_0, \Phi_1$ are two harmonic GBC maps from $V \mapsto \Theta$ and $W 
\mapsto \Theta$, respectively. 
The researchers in \cite{SH15} simply used the Rad\'o-Kneser-Choquet theorem 
to explain that  $\Phi_0$ is 
bijective from $V \mapsto \Theta$ and so is $\Phi_1$. 
Note that the classic Rad\'o-Kneser-Choquet theorem only explains the bijectivity of the unit ball
to a convex domain. See the wiki page of the  Rad\'o-Kneser-Choquet theorem on-line. 
See \cite{AN08} for a recent revisit of the   Rad\'o-Kneser-Choquet theorem. 
If one combines it with a Riemann map from a polygonal domain $V$ to the unit 
ball $B$, one can establish 
the bijection between an arbitrary polygon $V$ to a convex function $W$. Unfortunately, 
the  Rad\'o-Kneser-Choquet theorem does not hold for a higher dimensional setting other 
than $\mathbb{R}^2$. One can not use the  Rad\'o-Kneser-Choquet theorem 
for the bijectivity of  3D harmonic GBC maps.   
Although the researchers in \cite{SH15} have many nice examples in the  $\mathbb{R}^3$ setting, 
it is necessary to study the bijectivity of harmonic mappings in a general dimensional space. 
These motivate us to study the bijectivity of the harmonic map based on GBC functions from 
polygon $V$ to a convex polygon $W$. For simplicity, we consider only the situation in the 
Euclidean space $\mathbb{R}^2$ in this paper.   

\begin{figure}
\centering
\includegraphics[width= 0.8\textwidth]{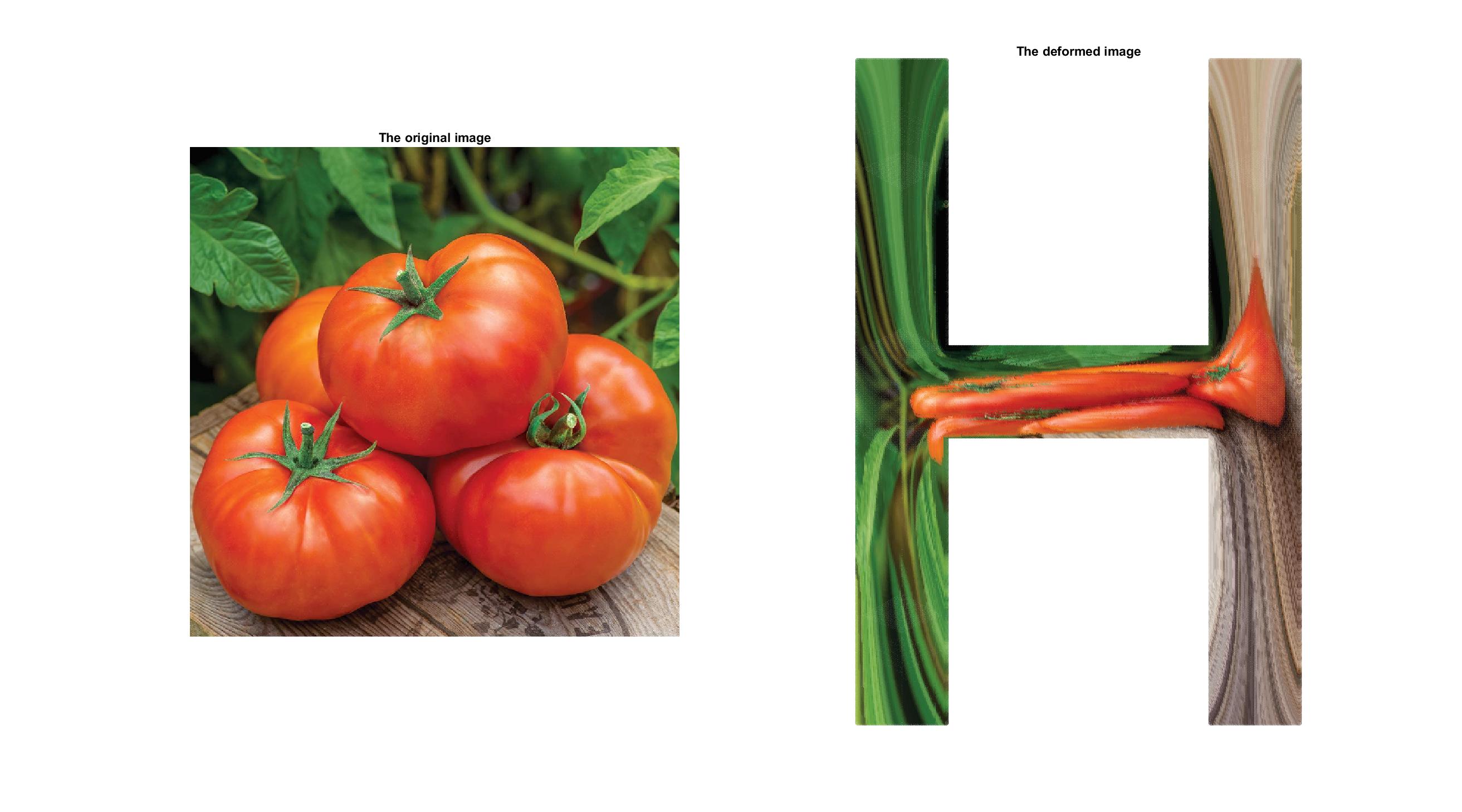}
    \caption{An  Example of   Image Deformation Based on Bijections from Polygons to Rectangle. 
Clearly the polygon on the right has more than 4 vertices. We added inactive vertices on 
the left domain to match the number of the vertices of the domain on the right.   \label{fig1}}
\end{figure}

\begin{figure}
\centering
\includegraphics[width= 0.8\textwidth]{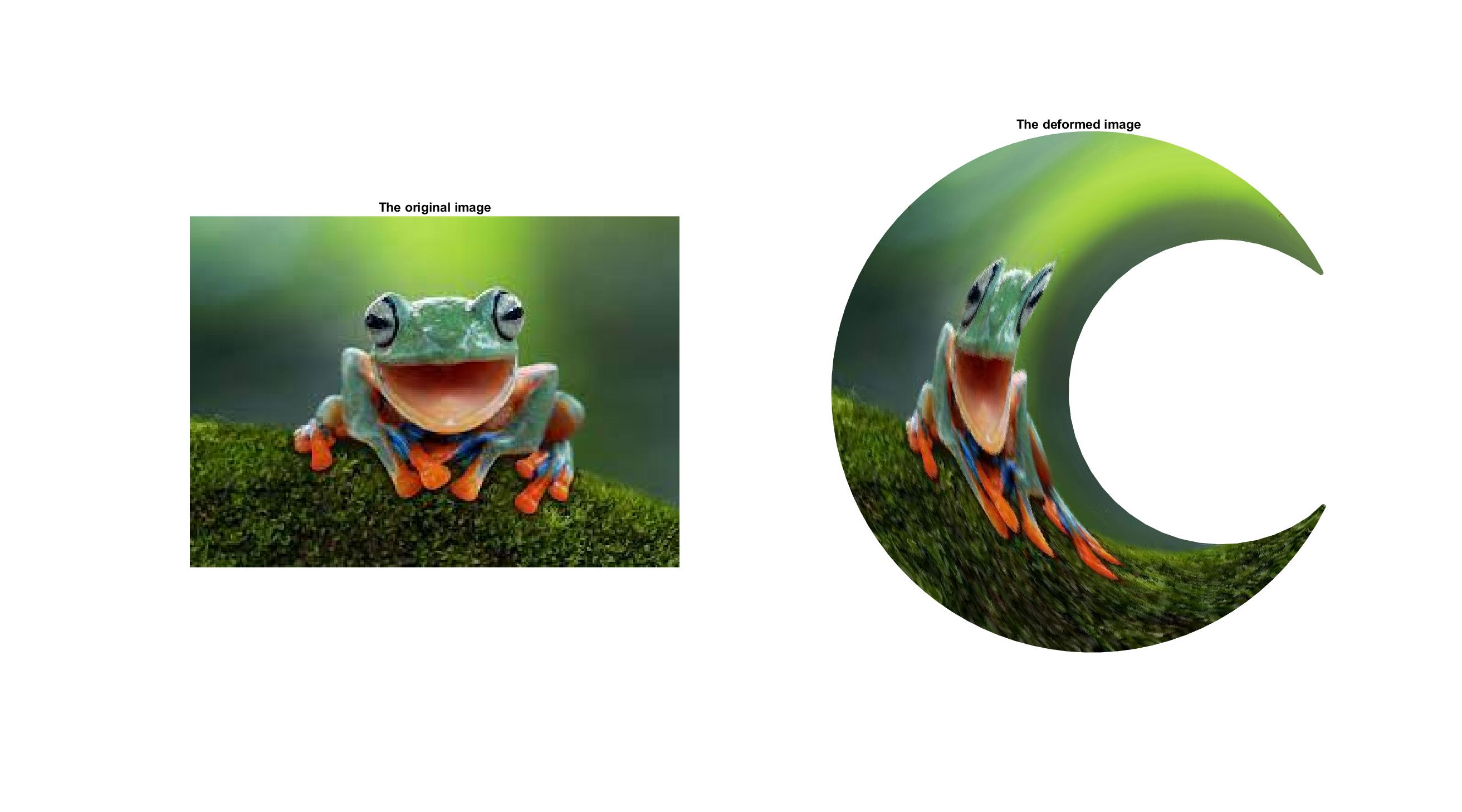}
    \caption{An Example of   Image Deformation Based on Bijections from Polygons to Rectangle. Note that the polygon on the right has more than 4 vertices.   \label{fig1b}}
\end{figure}

\begin{figure}
\centering
\includegraphics[width= 0.8\textwidth]{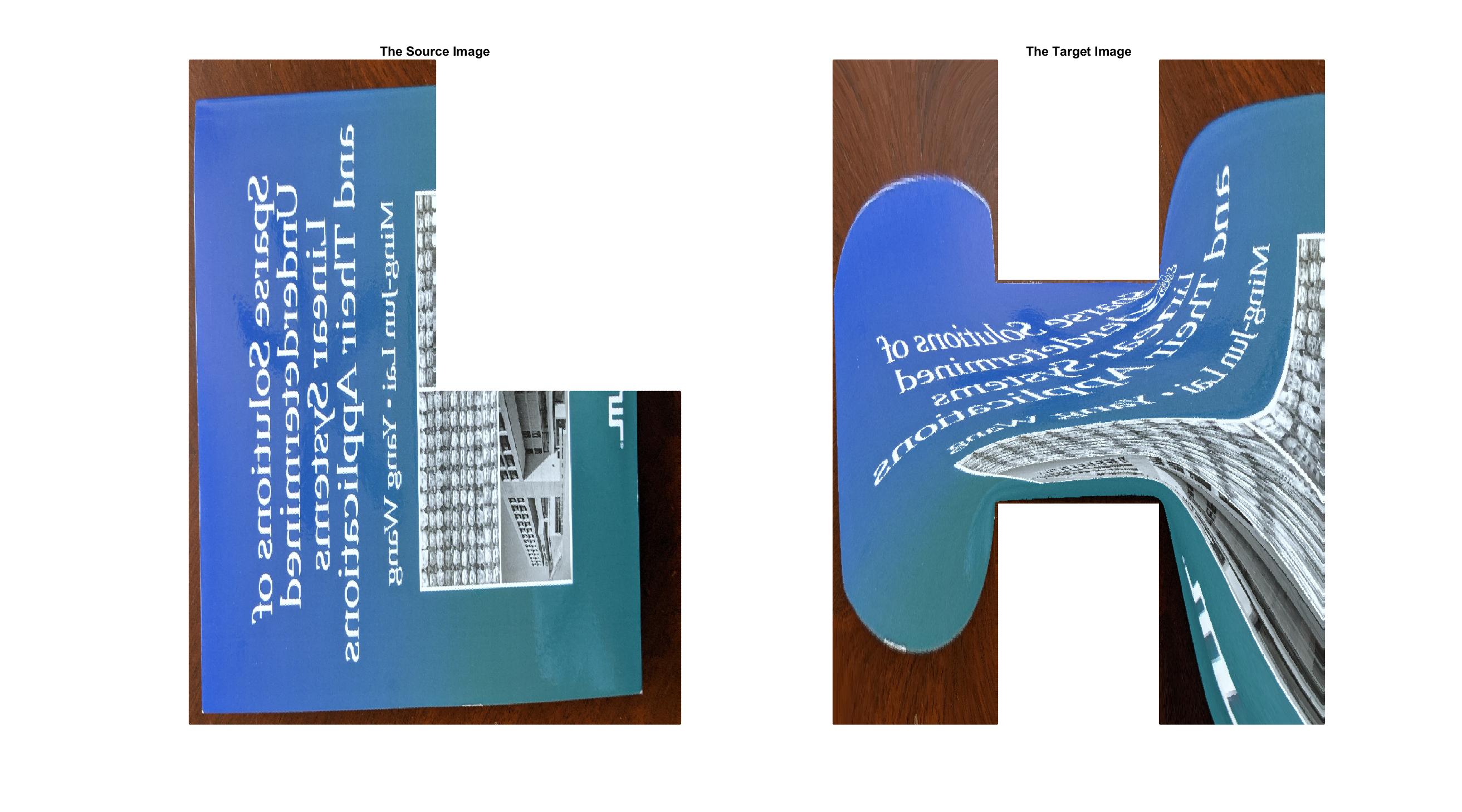}
\caption{An Example of Image Deformation from an L-shaped Polygon to the H-shaped polygon via a square domain based 
on the approach in \cite{SH15}\label{L2H}}
\end{figure}

The purpose of  this paper is to give a proof of the bijective property
of $\Phi_i$, $i=0, 1$ in the $\mathbb{R}^2$ setting. See Theorem~\ref{harbi}.   
In addition, we extend the study to  the bijectivty of a general harmonic map under some 
additional boundary condition. See Theorem~\ref{harbi2} in the next section.  
We leave the study of the nontrivial 
extension of the results in this paper to the multidimensional setting in a future publication.  
Finally, we point out that 
polygons $V$ and $W$ of  interest may not need to have the same number 
of  vertices.  One may add  inactive vertices to the boundary
of  $V$ or $W$ so that the number of  vertices of  $V$ and $W$ will be the same and 
the computation of  the bijection 
${\bf F}$  can be carried out based on the method in \cite{SH15}.  
 See three numerical examples in Figure~\ref{fig1}-- Figure~\ref{L2H} 
which are computed by a similar method to the one in \cite{SH15} 
 with a difference that we employed bivariate spline approximation of  
harmonic GBC functions discussed in  \cite{DFL20}.

From these images we can see how the points on the right-hand sides in Figure~\ref{fig1} and
Figure~\ref{fig1b} are mapped to the points on the left-hand sides by using the harmonic GBC
map.

In addition to the proof of Theorems~\ref{harbi} and \ref{harbi2}, we also make a few comments 
in the end of this paper. One example will be shown that a  
 harmonic GBC map from one nonconvex domain to another nonconvex domain can be bijective. 
Harmonic GBC functions can be also constructed over a polygonal domain with hole or holes. 
Does the harmonic GBC map based on these GBC functions have a bijective property?   
Besides the bijectivity, we will 
mention that the harmonic GBC map is also a diffeomorphism when $W$ is convex.  Although the 
harmonic GBC map does not have the property like a conforming map, the harmonic GBC map can be
easily computed since  harmonic GBC functions over all boundary vertices of $V$ can be done 
in parallel and once these functions are found, they can be used to form a bijection to any 
convex set $W$. A similar study of the bijectivity of harmonic GBCs in the 3D setting will be 
commented and difficulties will be pointed out. It seems the extension to the 3D setting is 
nontrivial.    

\section{Main Results and Their Proofs}
In this section,  let us start with harmonic functions over  simply connected domains in 
$\mathbb{R}^2$.  Let $V$ and $W$ be two simply connected domains. Typically, $V$ and $W$ 
are polygons.  
Let  $u: V \rightarrow \mathbb{R}$ is a harmonic function, that is, $\Delta u=0$ over 
$\mathring{V}$, where $\Delta = D_x^2+D_y^2$  is the standard Laplace operator. For two 
harmonic functions $u, v$ on $V$, let 
${\bf F}= (u(x,y), v(x,y))$ which will be called a harmonic map on  $V$ which maps $V$ to 
$\hbox{Im}({\bf F})$, the image of ${\bf F}$ over $V$.   

The harmonic functions have been studied for more than 100 years and many properties are 
known.  Let us mention a few properties without proofs which will be used in the following study.  
The harmonic function $u$ will not achieve its maximum nor minimum inside $V$ 
except for a constant function $u$. When $u$ is a harmonic, so are its derivatives. 
If a harmonic function $u$ is a constant over
a disk $D$ (open set) inside $V$, then $u$ is a constant over $V$.  Indeed, $\displaystyle 
\frac{\partial u}{\partial x}$ over $V$ will be zero over $D$ and hence $\displaystyle 
\frac{\partial u}{\partial x}$ is zero over $V$ as a harmonic function is like an 
analytic function. 

Note that when  $V, W \subset \mathbb{R}^2$ be two separate polygons,   
$\partial V$ and $\partial W$ are piecewise linear boundaries, respectively. 
By the properties of GBC functions 
(\ref{GBC}) and the definition of ${\bf F}$ (\ref{GBCmap}), 
one can easily check that ${\bf F}|_{\partial V}$ is a homeomorphism from 
$\partial V$ to $\partial W$.
We can also see that ${\bf F}|_{\partial V}$ preserve the orientation of the boundary 
vertices of $W$ as both of them are in the counter-clockwise direction. 

Our first main result in this paper is to establish the following 
\begin{theorem}\label{harbi}
Let $V, W\subset \mathbb{R}^2$ be two polygons.   
Let ${\bf F}=(u(x,y),v(x,y)): V\mapsto W$ 
be a harmonic GBC function which maps from  $V$ to $W$ according to (\ref{GBCmap}). 
If $W$ is convex, then ${\bf F}$ is a bijection from $V$ to $W$.
\end{theorem}

We remark that the above result can be generalized to the general harmonic map  
${\bf F}=(u,v): V \mapsto W$ under the assumptions that 
\begin{itemize}
    \item(1) ${\bf F}|_{\partial V}$ is an orientation preserving homeomorphism from the boundary $\partial V$ of $V$ to the boundary $\partial W$ of $W$.
    \item(2) ${\bf F} \in C^0(V) \cap C^1 (\mathring{V})$.
\end{itemize}

In this general setting, we have 
\begin{theorem}\label{harbi2}
Let $V, W\subset \mathbb{R}^2$ be two simply connected domains in $\mathbb{R}^2$.   
Suppose $u, v$ are  harmonic functions over $V$.  
Let ${\bf F}=(u(x,y),v(x,y)): V\mapsto W$ be a harmonic  map from  $V$ to $W$. 
Suppose that ${\bf F}$ satisfies (1) and (2) above.    
Then if $W$ is convex,  ${\bf F}$ is a bijection from $V$ to $W$.
\end{theorem}

Note that the map ${\bf F}$  is a bijection if ${\bf F}$  is surjective and also injective 
in the sense that ${\bf F}$  maps $V$ onto $W$ and the
mapping is one-to-one.  Before proving  Theorem~\ref{harbi}, we first establish the surjectivity.

\begin{lemma}
\label{onto}
As the setting in Theorem \ref{harbi2}, ${\bf F}$  maps $V$ onto $W$ if $W$ is convex. 
\end{lemma}
\begin{proof}
We will show that the image $\hbox{Im}({\bf F})$ of ${\bf F}$ over $V$ is $W$, i.e. 
$\hbox{Im}({\bf F})=W$.   

First of all, if ${\bf F}$ is a harmonic GBC, 
it is easy to see $\hbox{Im}({\bf F})  \subset W$ by the GBC properties (\ref{GBC}).  
In general, we claim $\hbox{Im}({\bf F})  \subset W$ for a general harmonic map ${\bf F}$. 
Otherwise,  there exists a point ${\bf p} = (u_0, v_0) \in 
\hbox{Im}({\bf F}) \setminus W$. 
Then there is a line separate $W$ and ${\bf p}$ because $W$ is convex. 
That is,  there exist real numbers $\alpha$ and $\beta$ such that 
$\alpha \cdot u_0 + \beta \cdot v_0 > 0$
, and $\alpha \cdot u' + \beta \cdot u' < 0$ for all $(u', v') \in W$. 
For $\bfx=(x,y)\in V$, let $\Phi(\bfx)=\alpha \cdot u(x,y) + \beta \cdot v(x,y)$ be an harmonic function over $V$ and 
${\bf F}^{-1}({\bf p}) \in \mathring{V}$. Then we have $\Phi({\bf F}^{-1}({\bf p})) >0$ 
and $\Phi(\bfx) <0$ for all ${\bf x} \in \partial V$,  this leads to a contradiction 
because it violate the maximum principle of the harmonic function $\Phi$. Hence, the claim 
follows.

On the other hand, we want to show that $W \subset \hbox{Im}({\bf F})$. 
To do so, we introduce the concept called 
the contractible curve. If a Jordan curve $\gamma$ 
inside $V$ can be continuously shrunk within $V$ to a point in $V$, then $\gamma$ is called
a contractible curve in $V$. We note that $V$ is simply connected, $\partial V$ is a 
contractible curve 
in $V$. Since ${\bf F}|_{\partial V}$ is a homeomorphism (continuously bijection) 
to $\partial W$, $\partial W$ must also be contractible in $\hbox{Im}({\bf F})$.

Now suppose $W \not\subset  \hbox{Im}({\bf F})$.  
Then there exists a point $ {\bf w} \in W \setminus \hbox{Im}({\bf F}) $.
Thus, $\partial W$ is contractible in $\hbox{Im}({\bf F}) \subset (W \setminus {\bf w})$. 
As $W\setminus {\bf w}$ is larger than $\hbox{Im}{\bf F}$,  $\partial W$ would be  
contractible in $(W \setminus {\bf w})$ which is obviously not possible as $W \setminus 
{\bf w}$ is like a punctured disk.        
Hence, $W \subset \hbox{Im}({\bf F}) $.
\end{proof}

Next we shall show that ${\bf F}$  is an injection.
Recall that ${\bf F}=(u,v)^\top$ is homeomorphism at $(x_0,y_0)\in V$ 
if its Jacobian matrix $J_0$ is not singular at the point $(x_0,y_0)\in \mathring{V}$. 
Suppose that ${\bf F}(x,y)$ has a nonzero Jacobian determinant at ${\bf x_0}=(x_0,y_0)$, that is,
\begin{equation}
\det( J_0 )= \det( [\nabla u(x_0, y_0); \nabla v(x_0, y_0)])= 
\det \begin{bmatrix} \frac{\partial }{\partial x}u(x_0,y_0) &\frac{\partial }{\partial y}u(x_0,y_0)\cr
\frac{\partial }{\partial x}v(x_0,y_0) & \frac{\partial }{\partial y}v(x_0,y_0)\end{bmatrix} 
\neq 0.  
\end{equation}

Since $\det (J_0) \neq 0$, $\left\| J_0 \cdot {\bf x}\right\| \ge c\left\| {\bf x} \right\|$ 
for some positive constant $c$ for all $\bfx$ in a neighborhood of $(x_0,y_0)$. 
The Taylor expansion at ${\bf x_0}=(x_0,y_0)$ gives us
\begin{eqnarray}
\left\| {\bf F}({\bf x})- {\bf F}({\bf x_0}) \right\|
&=&\left\| J_0 \cdot ({\bf x}-{\bf x_0}) \right\|+ o( \left\| {\bf x}-{\bf x_0} \right\|)\cr 
& \geq & c\left\| {\bf x}-{\bf x_0} \right\| + o( \left\| {\bf x}-{\bf x_0} \right\|).
\end{eqnarray}
Hence there exist a $r_0 >0$ such that ${\bf F}({\bf x}) \neq {\bf F}({\bf x_0})$ 
for all ${\bf x} 
\in B_{r_0}({\bf x_0})$. Now consider the closed ball $\overline{B}_{r_0}({\bf x_0}) \subset 
\mathring{V}$if $r_0>0$ small enough. Apply the same trick to each point on $\overline{B}_{r_0}({\bf x_0})$. We let 
$$r=\min \left\{ r_k \, | \, {\bf x_k} \in \overline{B}_{r_0}, {B}_{r_k}(\bfx_k) \subset 
\mathring{V}, {\bf F}({\bf x}) \neq {\bf F}({\bf x_k})\, for\, all\, {\bf x} \in B_{r_k}({\bf 
x_k}) \right\}.
$$
The compactness of $\overline{B_{r_0}(\bfx_0)}$ implies that 
$r>0$.  Then it is clear that ${\bf F}$ restricted on $B_{r/2}({\bf x_0})$ 
is a homeomorphism from $B_{r/2}({\bf x_0})$ to its image $\hbox{Im}({\bf F})$ over the ball.

The discussion above leads to the following
\begin{lemma}
\label{localinjection}
Suppose that the Jacobian determinant of ${\bf F}$ is not zero at every point of $V$. 
Then ${\bf F}$ is a local homeomorphism on $\mathring{V}$, 
i.e.. continuous, one-to-one and onto map over a neighborhood $N(x,y)\subset V$ 
for $(x,y)\in \mathring{V}$.
\end{lemma}

Unfortunately, the nonzero Jacobian determinant of
${\bf F}$ over $V$ does not guarantee the global injectivity as shown in the following example:
\begin{example}
Let ${\bf F}= (u,v)$ with $u(x,y)=  \exp(x/2) \cos(y \exp(-x))$ and $v(x,y)= \exp(x/2) 
\sin(y\exp(-x))$. Then the Jacobian determinant is equals to $1/2$ for all 
$(x,y)\in \mathbb{R}^2$. However, ${\bf F}(0,0)={\bf F}(0,2\pi)$. 
i.e. ${\bf F}$ is not injective globally.    
\end{example}

On the other hand, due to the GBC property, we can see that when $W$ is a triangle, the 
map ${\bf F}$ is an injection. Indeed, suppose that there are two points $(x_1,y_1)$ and
$(x_2,y_2)$ in $V$ which is also a triangle such that ${\bf F}(x_1,y_1)={\bf F}(x_2,y_2)$. 
Then we have
$$
\sum_{i=1}^3 (\phi_i(x_1,y_1)-\phi_i(x_2,y_2)) \bfw_i = 0.  
$$
Without loss of generality, we may assume that $\bfw_1=(0,0)$. As $W$ is a triangle with nonzero
area, it follows that $(\phi_i(x_1,y_1)-\phi_i(x_2,y_2))=0, i=2,3$. That is, $\phi_i(x_1,y_1)
= \phi_i(x_2,y_2), i=2, 3$. By the GBC property, $1= 
\phi_1(x_2,y_2)+\phi_2(x_2,y_2)+\phi_3(x_2,y_2)$ and hence, $\phi_1(x_1,y_1)=\phi_1(x_2,y_2)$. 
 Now 
$$
(x_1,y_1) =\sum_{i=1}^3 \phi_i(x_1,y_1)\bfv_i= \sum_{i=1}^3\phi_i(x_2,y_2) \bfv_i 
=(x_2,y_2).
$$ 
We see that $x_1=x_2, y_1=y_2$ and hence, ${\bf F}$ is injective. It is also known that when $W$ is not a convex 
polygon, the GBC map may not be an injection.  See a counter example in \cite{J13}. 

In the following, let us prove that the determinant of the Jacobian matrix of 
${\bf F}$ is nonzero if ${\bf F}$ is harmonic and $W$ is convex. To do so, 
we need more  properties of harmonic functions. 
Let us begin with the set of critical points of harmonic function.

\begin{definition}
Suppose that $u$ is harmonic on $\Omega \subset \mathbb{R}^2$.  
The critical set $C(u)$ of $u$ is the subset of $\Omega$ 
where the $\nabla u$ vanish. More precisely, 
$$C(u) \coloneqq \{ {\bf x} \in \Omega | \quad \nabla u ({\bf x})=0 \}. $$
\end{definition}

\begin{lemma}\label{isozero}
If $u$ is non-constant harmonic function, 
then $C(u)$ is an empty set or a collection of finitely many isolated points.
\end{lemma}

\begin{proof}
Let ${\bf i}=\sqrt{-1}$ and write $\bfz=x+{\bf i}y$.  
Denote by $u_x, u_y$ the partial derivatives of $u$.  
Note that $f({\bf z})=u_x-{\bf i} u_y$ is a holomorphic complex function as one can easily check 
that $f$  satisfied Cauchy-Riemann equations. A fundamental theory of complex analysis 
states that the set of zeros of an holomorphic function could only be empty or finitely many 
isolated point(s). Since $\nabla u (x,y) =0$ if and only if $f({\bf z})=0$, 
we have $C(u)= \{ {\bf x} \in \Omega| f({\bf z})=0 \}$ is empty or contains finitely many
isolated point(s).
\end{proof}

\begin{remark}
There are few different ways to prove Lemma \ref{isozero}. 
For example, a good way to see is that the 
critical set $C(u)$ of a harmonic function $u: \mathbb{R}^n \rightarrow \mathbb{R}$ have 
Hausdorff dimension at most $n-2$ as discussed in \cite{HHHN}. 
So when $n=2$, $C(u)$ is empty or isolated point(s) which are only cases when a set is
of  Hausdorff dimension $0$.
\end{remark}

Next we discuss the level curve(s) of harmonic function $u$. Fix $\bfa\in \Omega$. 
Let $\Gamma_u({\bf a}) := \{ {\bf x} \in \Omega | u({\bf x}) = u({\bf a}) \}$ 
be the level curve set of $u$ at $u({\bf a})$. 
The set will have the following properties. 
\begin{lemma}
\label{Property2}
Let $u: \Omega \rightarrow \mathbb{R}$ be a non-constant harmonic function. 
If ${\bf a} \in C(u) $, 
then $\Gamma_u(\bf a)$ are curves intersect at ${\bf a}$. Moreover, the angle between two neighboring curves is $ \pi /n$, where $n$ is the number of curves
 in $\Gamma_u({\bf a})$. Also, any of these curves is  not a loop.
\end{lemma}
\begin{proof}
WLOG, let ${\bf a}=0$. Consider a function $f({\bf z}) =u_x-{\bf i} u_y$. 
Then $f(0)=0$ and $f$  is 
holomorphic and hence, $f$  has a finite multiplicity, say $n$ at $0$. 
Write $f=nz^{n-1}+O(z^n)$ be the  
Taylor expansion of $f$ at $0$. 
So $u=Re(z^n)+O(z^{n+1})$. Note that $\Gamma_{z^n}(0)$ is $n$ lines 
intersects at $0$ with angle $\frac{2 \pi}{n}$, 
so  $\Gamma_{u}(0)$ are curves as described above. 

Next if there is a level curve which has a loop, say it encloses a nonempty open set  
$D\subset V$. Because $u$ is harmonic on $V$, $u$ is harmonic on $D$ and hence, $u$ must be a
constant on $D$. That is, the critical point set of $u$ contains a nonempty open set $D$. 
which contradicts to Lemma \ref{isozero}. 
$u$ can have only isolated critical points. 
Hence, any curve in the level set is not a loop.
\end{proof}

Next we need the set of extremes.  
\begin{definition}
Let $ Max_u := \{ {\bfx} \in \partial \Omega | u({\bf x}) \geq u({\bf y}) 
\forall {\bf y}$ in a   
neighborhood of $\bfx \}$ be the set of 
all local maximal locations. Similarly, let $ Min_u := \{ \bfx \in \partial \Omega | u(\bfx) 
\leq u(\bfy) \forall \bfy$ in a neighborhood of $\bfx \}$ be the set of  
local minimal locations. 
We denote by $\# Max_u$ the number of connected components in the  set $Max_u$. 
Similarly,  $\# Min_u$ is the number of connected components in the set $Min_u$.   
\end{definition}

\begin{theorem}\label{cu}
Suppose that $u: \Omega \rightarrow \mathbb{R}^2$ is a harmonic function. If $\# Max_u = \# 
Min_u =1$, then $C(u)=\emptyset$.
\end{theorem}
\begin{proof}
If there exist any ${\bf a} \in C(u)$, by Lemma~\ref{Property2}, $\Gamma_u(\bf a)$ are curves 
intersect at $\bf a$. Since ${\bf a} \in \Omega$, these curves will intersect with $\partial 
\Omega$ $2n$ times, where $n$ is number of level curves in $\Gamma_u(\bf a)$ and $n \geq 2$. 
That is, the level curves cut the boundary $\partial \Omega$ into 2n pieces with $n\ge 2$.   
There must have some components of $Max_u$ or $Min_u$ between two intersecting points. 
So $\# Max_u \geq 2$ or $\# Min_u\ge 2$ which is a contradiction.
\end{proof}

Furthermore, we need more definitions and one more useful lemma. 
Define a continuous function $f_{\alpha}=\cos (\alpha)u + \sin (\alpha)v$ over $V$
for each $\alpha \in [0, \pi)$.   Let 
$$S_\alpha:=\{(u,v) \in W | \cos (\alpha)u + \sin (\alpha)v = \max_{(u',v') \in W} 
\cos (\alpha)u' + \sin (\alpha)v' \}$$
and
$$s_\alpha:=\{(u,v) \in W | \cos (\alpha)u + \sin (\alpha)v = \min_{(u',v') \in W} 
\cos (\alpha)u' + \sin (\alpha)v' \}.$$

\begin{lemma}\label{falpha}
 Let ${\bf F}$ , $V$ and $W$ be the same in  Theorem \ref{harbi2}. 
Then $\# Max_u = \# Min_u =1$ for any $\alpha \in [0, \pi)$.   
\end{lemma}
\begin{proof}
We denote by $\# S$ the number of connected components in the  set $S$ in the  proof. 
Fixed an $\alpha$ in $[0, \pi]$. 
Note that $f_{\alpha}$ is the inner product $(\cos \alpha ,\sin \alpha)\cdot {\bf F} $. 
See an illustration of $f_\alpha$ in Figure~\ref{ex1}. 

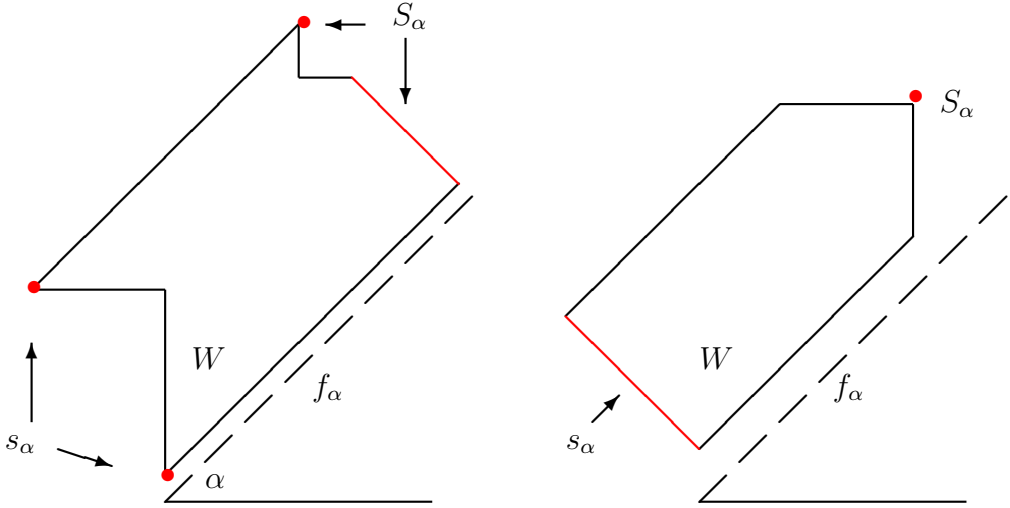
\begin{figure}
\begin{picture}(0, 200)
\thicklines
\put(50,80){\line(1,1){100}}
\put(100,10){\line(0,1){70}}
\put(50,80){\line(1,0){50}}
\put(110,50){$W$}
\put(100,10){\line(1,1){110}}
\put(98,7){\textcolor{red}{$\bullet$}}
\put(48,78){\textcolor{red}{$\bullet$}}
\put(40,20){$s_\alpha$}
\put(50,30){\vector(0,1){30}}
\put(60,20){\vector(3,-1){20}}
\put(150,180){\line(0,-1){20}}
\put(149,178){\textcolor{red}{$\bullet$}}
\put(150,160){\line(1,0){20}}
\put(170,160){\textcolor{red}{\line(1,-1){40}}}
\put(185,180){$S_\alpha$}
\put(175,180){\vector(-1,0){15}}
\put(190,175){\vector(0,-1){25}}
\put(100,0){\line(1,0){100}}
\put(115,5){$\alpha$}
\multiput(100,0)(15,15){8}{\line(1,1){10}}
\put(155,40){$f_\alpha$}

\put(300,0){\line(1,0){100}}
\multiput(300,0)(15,15){8}{\line(1,1){10}}
\put(350,40){$f_\alpha$}
\put(300,20){\line(1,1){80}}
\put(250,70){\line(1,1){80}}
\put(300,50){$W$}
\put(300,20){\textcolor{red}{\line(-1,1){50}}}
\put(380,100){\line(0,1){50}}
\put(330,150){\line(1,0){50}}
\put(378,150){\textcolor{red}{$\bullet$}}
\put(250,20){$s_\alpha$}
\put(260,30){\vector(1,1){10}}
\put(390,147){$S_\alpha$}
\end{picture}
\caption{An illustration of the sets ${S_\alpha}$ and
$s_\alpha$ which have multiple components  
when $W$ is  a nonconvex (left) and the sets 
$S_\alpha$ and
$s_\alpha$ when $W$ is a convex domain 
(right) \label{ex1}}
\end{figure}

Since $W$ is convex, both $S_\alpha$ and $s_\alpha$ are a set 
consisting of one  connected component. 
See Figure~\ref{ex1} for the cases when $W$ is convex and is not convex. When $W$ is convex, 
$S_\alpha$ contains only one component of the image ${\bf F}$. 
When $W$ is not convex, $S_\alpha$ may
have two or more separate components as in the left of Figure~\ref{ex1} where $\alpha=\pi/4$.  
Because  $f_\alpha$ is harmonic, $f_\alpha$ 
can only achieve its maximum on the boundary. So $S_\alpha \subset \partial W$. 
Because ${\bf F}$  is a continuous  bijection on the boundary $\partial V$ to $\partial W$, 
${\bf F}^{-1} (S_\alpha)$ is a set of one component and hence $\# Max_u =\# {\bf F}^{-1} (S_\alpha) =1$.  
Similarly, $\# Min_u=\# {\bf F}^{-1} (s_\alpha)=1$. 
\end{proof}

Note that even though the non-zero determinant of Jacobian implies local injection, we cannot 
extend it to $\partial V$ because ${\bf F}$ is not necessary be differentiable on $\partial V$ 
as $V$ has corners. However, we can prove the following 
\begin{lemma}
\label{finite}
Suppose that ${\bf F}: K \rightarrow T$, where $K \subset V$ is a compact set and $T \subset W$ 
is connected. Also, ${\bf F} \in C^0(K)$ is a local homeomorphism. 
The pre-image set ${\bf F}^{-1}(u,v)$ is finite for all $(u,v) \in \hbox{Im}({\bf F})$. 
\end{lemma}
\begin{proof}
Suppose that there is a point $(u,v)\in T$ such that the set ${\bf F}^{-1}(u,v)\subset K$ has an 
infinitely many points. Then due to the compactness of $K$, there is a subset of points which 
converge to a cluster
point $(x,y)\in K$. Then ${\bf F}$ is no longer a local homeomorphism at $(x,y)$ which 
contradicts to the result in Lemma~\ref{localinjection}.   
\end{proof}

Due to the fact that the pre-image set ${\bf F}^{-1}(u,v)$ consists of 
finitely many isolated points,  we can further show 
\begin{lemma}\label{l2}
Under the same assumptions in Lemma~\ref{finite}, 
 the number of all isolated points in the pre-image set ${\bf F}^{-1}(u,v)$ is a constant 
overall $(u,v) \in 
\mathring{W}$.
\end{lemma}
\begin{proof}
Suppose $T \subset \mathring{W}$ is a connected compact set. Since ${\bf F}$ is continuous, 
$K:={\bf F}^{-1}(T)$ is a compact set. ${\bf F}|_K$ is a local homeomorphism because 
$K \subset \mathring{V}$. 
We first show that $T_i := \{ (u,v) \in T | \# {\bf F}^{-1}(u,v)=i \}$ 
is an open set for any fixed integer $i\ge 1$, where $\# {\bf F}^{-1}(u,v)$ stands for the 
cardinality of the pre-image set $\# {\bf F}^{-1}(u,v)$ of $(u,v)\in T$. 

Let $(u,v) \in T_i$. Consider ${\bf F}^{-1}(u,v)=\{ \bfx_1, \bfx_2, \cdots, \bfx_i \}$, 
where $\bfx_j$ are isolated points in $\mathring{V}$. 
Since ${\bf F}$ is a local homeomorphism, we can choose disjoint neighborhood $N_j$ of $\bfx_j$, 
$j=1,2, \cdots, i$, where ${\bf F}$ is a homeomorphism from $N_j$ to the image set 
${\bf F}(N_j)$. 
Define $R=V \setminus (\cup_{j=1}^i N_j) $. Then $R$ is a close set as $V$ is compact and hence, 
${\bf F}(R)$ is a 
close set. Now we define an open set $M=(\cap_{j=1}^i {\bf F}(N_j)) \setminus {\bf F}(R)$ on 
$T$. Obviously, 
$M$ is a neighborhood of $(u,v)$, and ${\bf F}^{-1}(M)$ is a union of disjoint open sets belongs 
to $\cup_{j=1}^i N_j$. 
By the definition of $N_j$, ${\bf F}^{-1}(M) \cap N_j$ is a  homemorphism which means the 
cardinality 
$\# {\bf F}^{-1}((u',v'))=i$ for all $(u',v') \in M$, so $M \subset T_i$ and hence $T_i$ 
is an open set. 

Now since $T_i$ is open for any fixed $i$, it follows that $T={\bigcup}_i T_i$. 
As $T$ is connected, it cannot be a disjoint union of more than $2$ open set. This implies that
$T= T_i$ for some $i \geq 1$. Finally, because this is true for all connected set $T \subset 
\mathring{W}$, we can conclude that the number of the pre-image set ${\bf F}^{-1}(u,v)$ 
is a constant overall $(u,v) \in \mathring{W}$.
\end{proof}

\begin{figure}
\centering
\includegraphics[width= 0.7\textwidth]{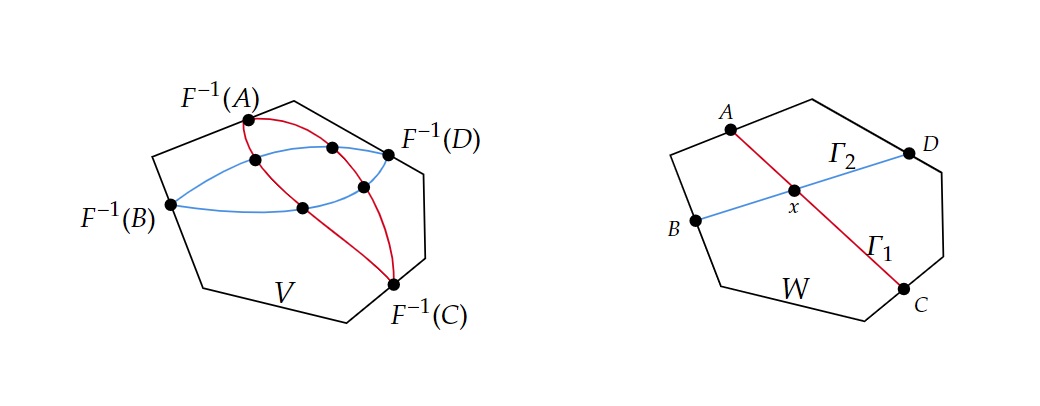}
\caption{An example when the cardinality of the pre-image set is 2 \label{TW}}
\end{figure}

We now show that $i=1$.  
Suppose that $i \ge 2$ on $\mathring{W}$. 
We pick 4 points $A$, $B$, $C$, and $D$ in the counterclockwise 
order as shown on the right of Figure~\ref{TW}. 
Then we pick a curve called $\Gamma_1$ connected $A$ and $C$, and a curve $\Gamma_2$ 
connected $B$ and $D$ over $W$, where 2 curves intersect at $x=(u,v)$ on the right of
Figure~\ref{TW}. Note that ${\bf F}$ is a 
homeomorphism on boundary, 
but not on the interior. We can see that 
${\bf F}^{-1}(\Gamma_1)$ has $i$ curves jointing ${\bf F}^{-1}(A)$ and 
${\bf F}^{-1}(C)$ and do not intersects over $\mathring{V}$. Similarly, 
${\bf F}^{-1}(\Gamma_2)$ has $i$ curves jointing  
${\bf F}^{-1}(B)$ and ${\bf F}^{-1}(D)$ and do not intersects over $\mathring{V}$. 
By the fact that the intersections of ${\bf F}^{-1}(\Gamma_1)$ and ${\bf F}^{-1}(\Gamma_2)$ 
is ${\bf F}^{-1}(u,v)$. 
We conclude that the cardinality $\# {\bf F}^{-1}(u,v)=i^2>i$, which is a violation of 
Lemma~\ref{l2} unless $i=1$. Figure~\ref{TW} shows the case when $i=2$.

Let us summarize the discussion above in the following 
\begin{theorem}
\label{thm1}
Suppose $V$, $W$ are two polygons in $\mathbb{R}^2$. 
Let ${\bf F}: V \rightarrow W$ be a harmonic GBC map. 
Then ${\bf F}$ is an injection from $V$ to $W$.
\end{theorem}
\begin{proof}
As we have seen that $i=1$, i.e. the number of isolated points in the pre-image
${\bf F}^{-1}(u,v)$  is one for each $(u,v)\in \mathring{W}$, the map ${\bf F}$ is injective.   
\end{proof}

Combine with Lemma \ref{localinjection} we have the following 
\begin{theorem}
\label{thm1b}
Suppose $V$, $W$ are two simply connected domains in $\mathbb{R}^2$. 
Let ${\bf F}: V \rightarrow W$ be a harmonic map satisfying 
${\bf F} \in C^1(\mathring{V}) \cap C^0(V)$. 
Suppose that ${\bf F}|_{\partial V}= \partial W$ is an orientation preserved homeomorphism and 
determinant of Jacobian is non-vanish on $\mathring{V}$. 
Then ${\bf F}$ is an injection from $V$ to $W$.
\end{theorem}

Now we ready to prove Theorem \ref{harbi2}.
\begin{proof}[Proof of Theorem \ref{harbi2}]
We write $(u,v)= {\bf F}(x,y)$ and use $f_\alpha$ for $\alpha\in [0,\pi)$.   
By Lemma \ref{falpha}, and Theorem \ref{cu}, 
we have $C(f_{\alpha})=\emptyset$. This implies the vector 
$$
\nabla f_\alpha = \nabla {\bf F}(x,y)\begin{bmatrix} \cos(\alpha)\cr \sin(\alpha)\end{bmatrix} \not=0
$$ 
for all $(x,y) \in \mathring{V}$ and for any $\alpha \in [0, \pi)$. 
So $\nabla u$ and $\nabla v$ are linearly independent in $V$. 
Hence, $\det(\nabla {\bf F})=\det([\nabla u; \nabla v]) \neq 0$ 
for all $(x,y) \in \mathring{V}$, which implies that ${\bf F}$  is a local homeomorphism 
on $\mathring{V}$ by Lemma~\ref{localinjection}. 

Hence, ${\bf F}$ is an injection by  Theorem~\ref{thm1} if ${\bf F}$ is a harmonic 
GBC map and by Theorem~\ref{thm1b} if ${\bf F}$ is a general harmonic map.  
Together with Lemma~\ref{onto}, we complete the proof.  
\end{proof}

\section{Remarks}
We have the following remarks in order.
\begin{itemize}
\item 1. Note that the convexity condition in Theorem~\ref{harbi} is a sufficient condition. 
It is interesting to see how the map ${\bf F}$ from $V$ to $W$ behaviors when 
$W$ may not be a convex domain.  We use our bivariate spline harmonic map ${\bf F}$ 
(cf.  \cite{DFL20}) to explore the possibility if the map ${\bf F}$ is a bijection. 
It is easy to find many examples that ${\bf F}$ does not map into $W$ if $W$ is 
not a convex polygon. 
However we found a few examples of polygons which can be mapped into an L-shaped domain
bijectively using our ${\bf F}$ in the numerical sense. 
See Figure~\ref{fig2} for such an example. We show a J-shaped
domain with labeled vertices in $J$ and in $L$.   
We compute the determinants of Jacobi matrix numerically at more than
25,000 points and only found 3 points on the boundary of the $L$-shaped domain where 
the determinant of the Jacobi matrix is less than $1e-5$. These three
points locate on the boundary of the L-shaped domain as shown in the red points 
in Figure~\ref{fig2}.  The contour lines on the right 
of Figure~\ref{fig2} show that the mapping is surjective. 
However we are not able to show such a map is indeed a bijection theoretically. 
\begin{figure}[h]
\centering
\includegraphics[width= 0.4\textwidth]{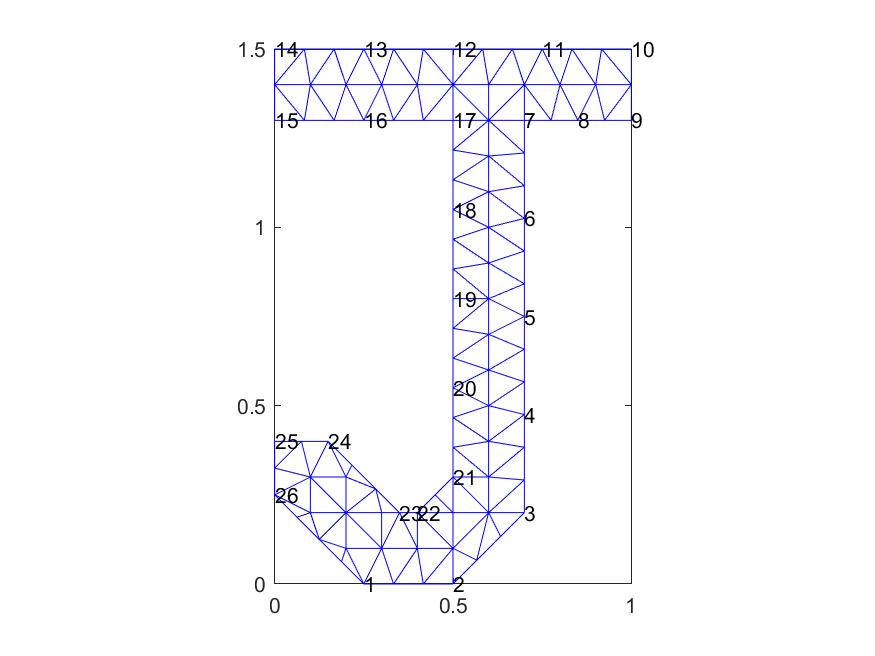}
\includegraphics[width= 0.4\textwidth]{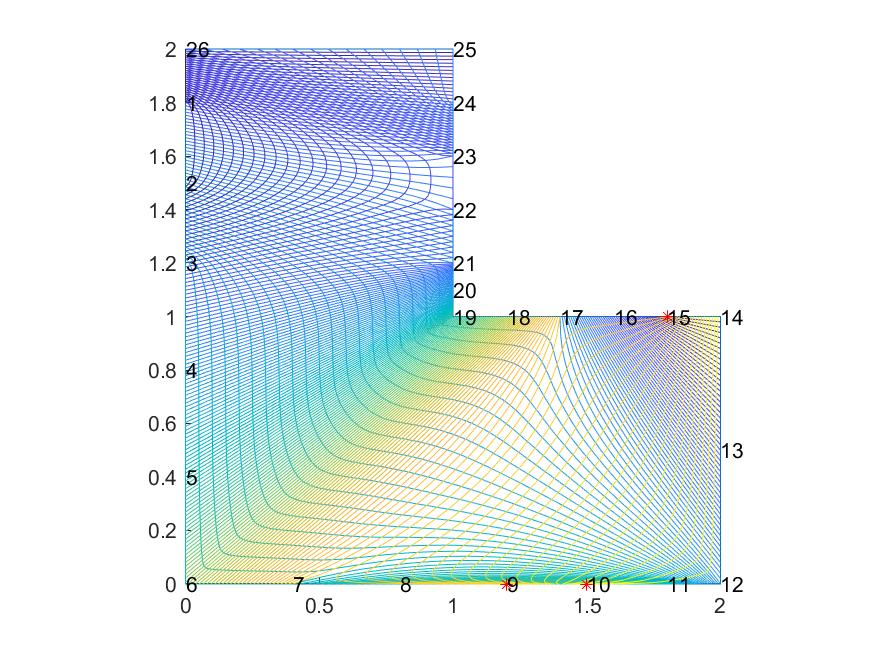}
\caption{Contour lines of ${\bf F}$ from a J-shaped domain to an L-shaped domain \label{fig2}}
\end{figure}

In fact, the mapping between the boundary of $V$ and the boundary $W$ is important. 
Many different arrangements of the boundary map between $V$ and $W$ lead to the 
harmonic GBCC map ${\bf F}$ which is not bijective in our numerical experiment. 
We leave the question how to find a good arrangement of the boundaries of $V$ and $W$ so  
that the map is bijective as an open problem for further study. 
 
\item 2. As harmonic GBC functions are infinitely many differentiable, we can conclude that 
${\bf F}$ defined in (\ref{GBCmap}) is not only a bijection, it also a diffeomorphism 
from $\mathring{V}$ to $\mathring{W}$.   
It is interesting to compare our results, e.g. Theorem~\ref{harbi} to the main result, i.e. Theorem 1.3 
in \cite{AN08}. For convenience, let us start the main result below. 
Let $B := \{(x, y) \in \mathbb{R}^2: x^2 + y^2 < 1\}$ denote the unit disk. 
Recall the following theorem from \cite{AN08}:
\begin{theorem}[G. Alexssandrini  and V. Nesi, 2008\cite{AN08}] 
Let $\Phi: \partial B \mapsto \gamma\subset \mathbb{R}^2$  
be an orientation preserving diffeomorphism
of class $C^1$ onto a simple closed curve $\gamma$. Let $D$ be the bounded domain such that
$\partial D = \gamma$. Let ${\bf U} \in C^2(\mathring{B}; \mathbb{R}^2)\cap  C^1(\overline{B}; 
\mathbb{R}^2)$ 
be the solution to the Dirichlet boundary problem of Laplace equation:  
\begin{eqnarray*}
\Delta {\bf U}(x,y) &=& 0, \quad (x,y)\in B \cr
{\bf U}(x,y) &=& \Phi, \quad (x, y)\in \partial B.  
\end{eqnarray*}
Then the mapping ${\bf U}$ is a diffeomorphism of $B$ onto $D$ if and only if
\begin{equation}
\label{keycondition}
 \det (\nabla {\bf U}) > 0 \hbox{ everywhere on } \partial B.
\end{equation}
\end{theorem}

Our result uses a convex domain $W$ instead of the unit $B$. 
Also, our ${\bf F}$ maps from a polygon $V$ to $W$ instead of ${\bf U}$
from $B$ to $D$ in \cite{AN08}. Our ${\bf F}$ is an orientation preserving homeomorphism 
over the piecewise linear boundary $\partial
V$ while the theorem above requires the diffeomorphism $\Phi$. 
We do not have a condition similar to (\ref{keycondition}). One significance of the 
harmonic GBC map ${\bf F}$ is that the diffeomorphism from $V$ to $W$ can be easily constructed using the GBC functions. 


\item 3. It is also interesting to compare the harmonic GBC map with the Riemann map from 
any polygon to the disk. The computation of such a Riemann map can be found in 
book \cite{GY07}. Both maps are diffeomorphism 
from $V$ to $W$ when $W$ is a disk. It is well-known that the Riemann map is a conforming map 
which preserves the angle and relative locations of the domain $V$ as shown in the images inside  
\cite{GY07}. However,  a harmonic GBC map  gives a distorted image called image warping 
as shown in Figure~\ref{L2H} which may be used for creative artworks such as cartoon images, 
e.g. in \cite{JMDGS06}.  One advantage
of the harmonic GBC map is the convenience of the construction: once GBC functions over a 
polygon $V$ are computed (which can be done in parallel, e.g. GPU), 
they can be used to map bijectively to any convex polygon $W$.    
In addition, the map between the boundary of $V$ and the boundary of $W$ can be chosen 
by the user which is a convenience as the Riemann map does not have this freedom.    

\item 4. Although the harmonic GBC functions are defined over a polygon, 
they are possible to be defined
over a polygonal domain with a hole or holes. We refer to \cite{HF06} for mean value 
coordinates over domains with holes.  
Indeed, suppose a polygonal domain $P$ with a hole
has $n$ outer boundary vertices $\bfv_i, i=1, \cdots, n$ and $m$ vertices  
$\bfh_j, j=1, \cdots, m$ on the inner boundary of the hole. 
One simply defines $\phi_1, \cdots, \phi_n$ are nonnegative functions 
which satisfy the standard GBC boundary 
conditions with additional property: 
$\phi_i$'s are zero at the boundary of the hole. 
Also, $\phi_{n+i}, \cdots, \phi_{n+m}$ are nonnegative functions  
which are zero on the outer boundary of $P$ and satisfy the 
standard GBC boundary conditions over the boundary 
of the hole. All these functions  $\phi_1, \cdots, \phi_n$ and $\phi_{n+i}, \cdots, 
\phi_{n+m}$ satisfy the standard GBC properties inside $P$:
$$
\sum_{i=1}^{n+m} \phi_i  = 1  
\hbox{ and }  \sum_{i=1}^n \phi_i \bfv_i + \sum_{j=1}^m \phi_{n+j}\bfh_j =(x,y),
$$ 
where $\bfv_i, i=1, \cdots, n$ are the vertices of 
the outer boundary of $P$ and $\bfh_j, i=1, \cdots, m$ are 
the vertices of the inner boundary of $P$.  
See Figs.~\ref{gbcholes} and ~\ref{gbcholes2} for  contour plots 
of some harmonic GBC functions ($\phi_{17}$ is now shown) computed based on 
the bivariate spline solution to the corresponding 
Laplace equation over the polygonal domain $P$ with a hole with an appropriate boundary 
condition (cf. \cite{ALW06}). 
It is interesting to know the bijectivity of the harmonic map based on these GBC functions. 
We leave the question to an open problem. 
\begin{figure}[h]
\centering
\includegraphics[width= 0.2\textwidth]{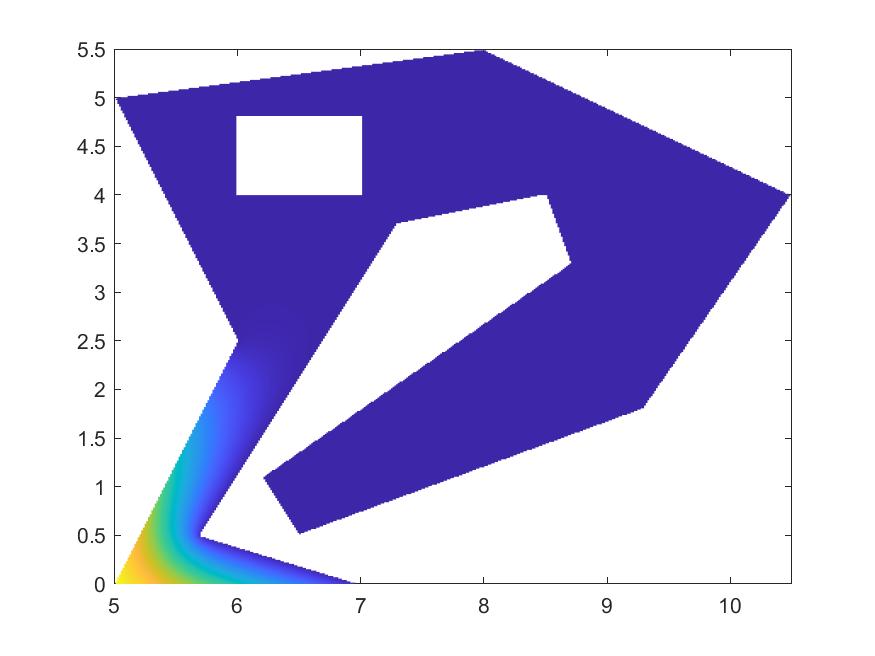}
\includegraphics[width= 0.2\textwidth]{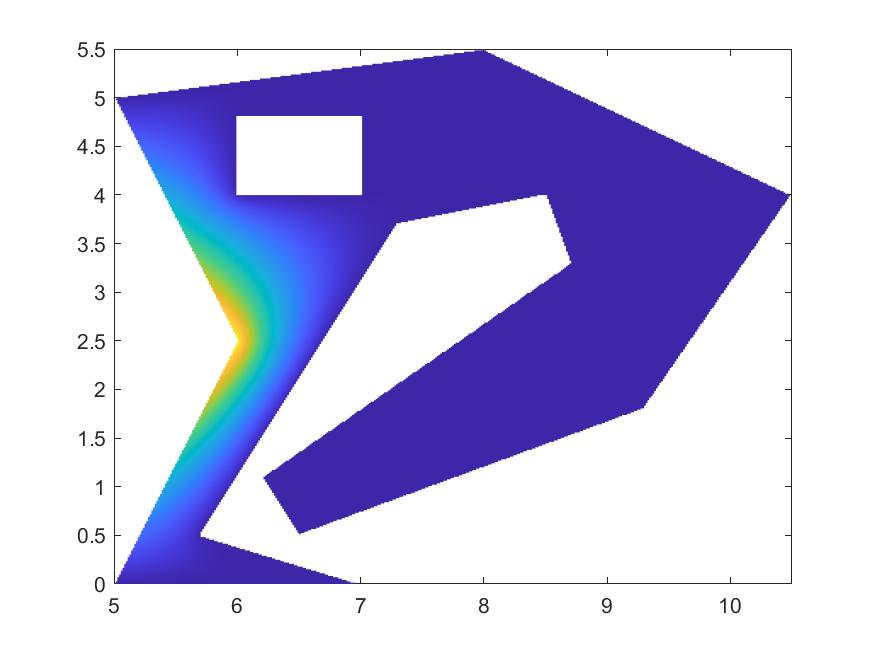}
\includegraphics[width= 0.2\textwidth]{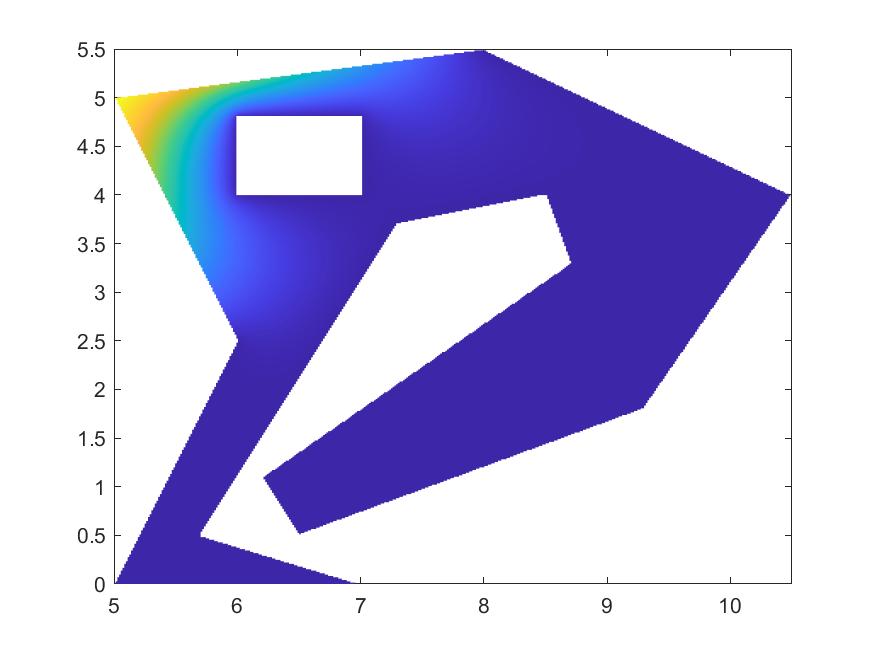}
\includegraphics[width= 0.2\textwidth]{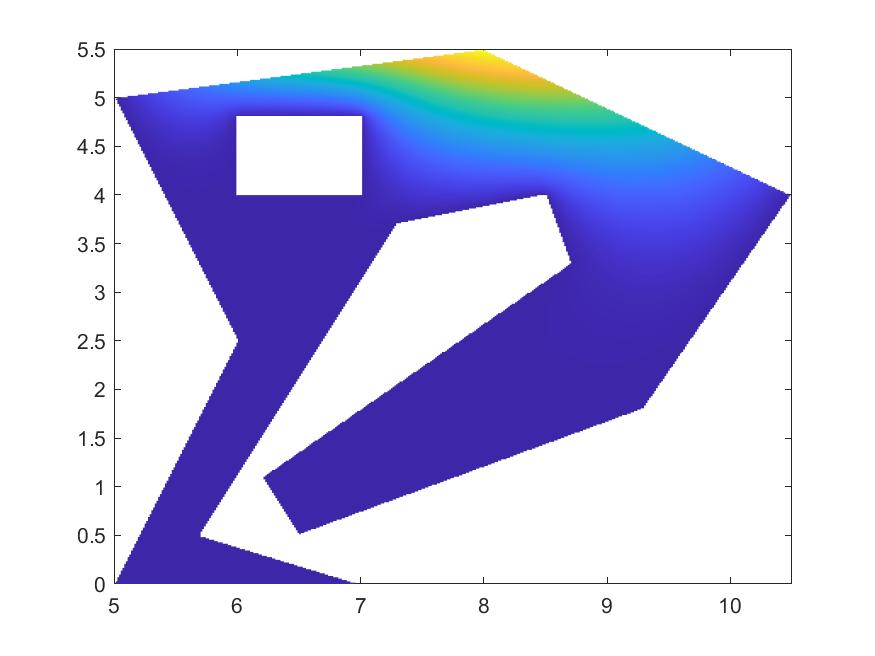}
\includegraphics[width= 0.2\textwidth]{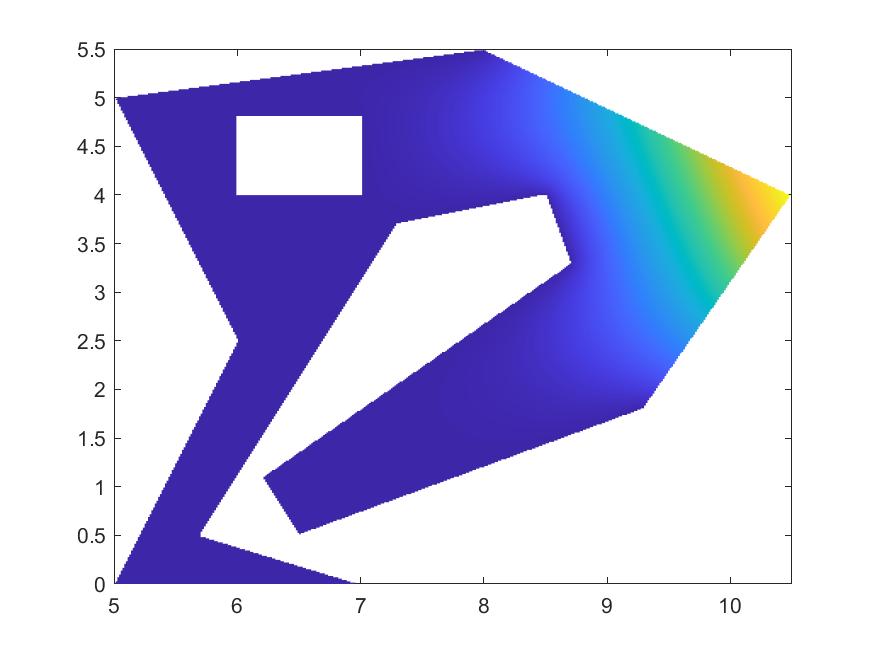}
\includegraphics[width= 0.2\textwidth]{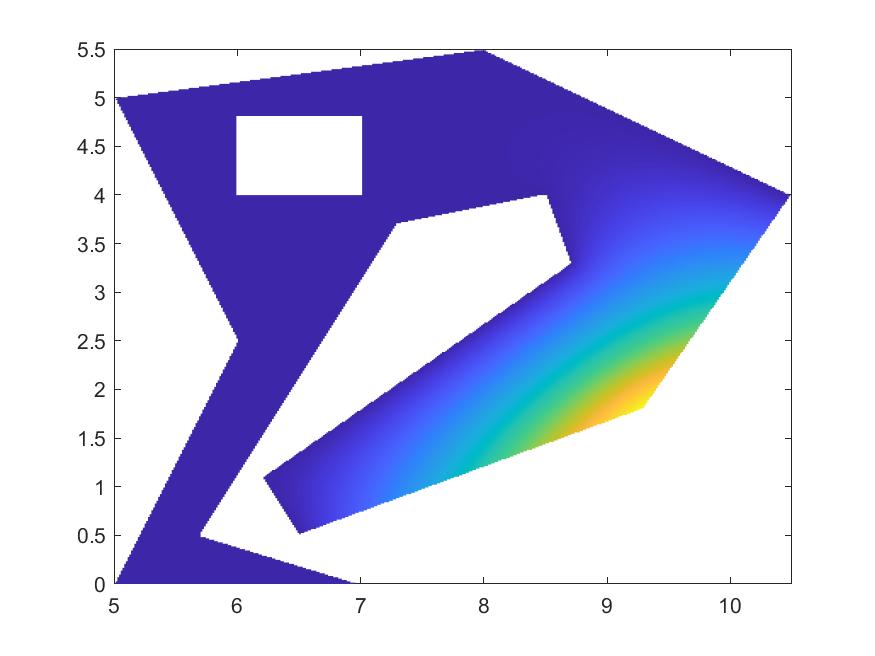}
\includegraphics[width= 0.2\textwidth]{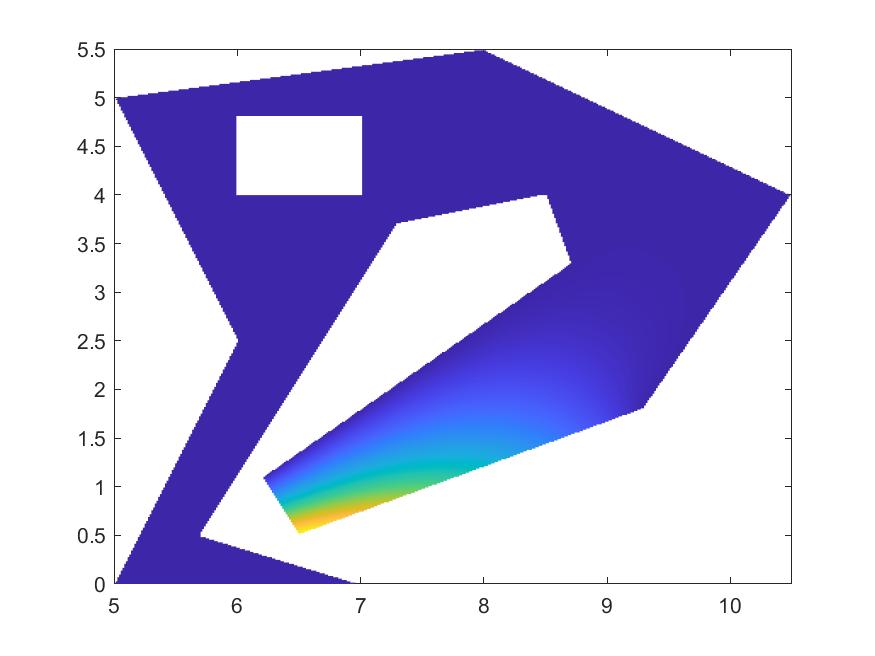}
\includegraphics[width= 0.2\textwidth]{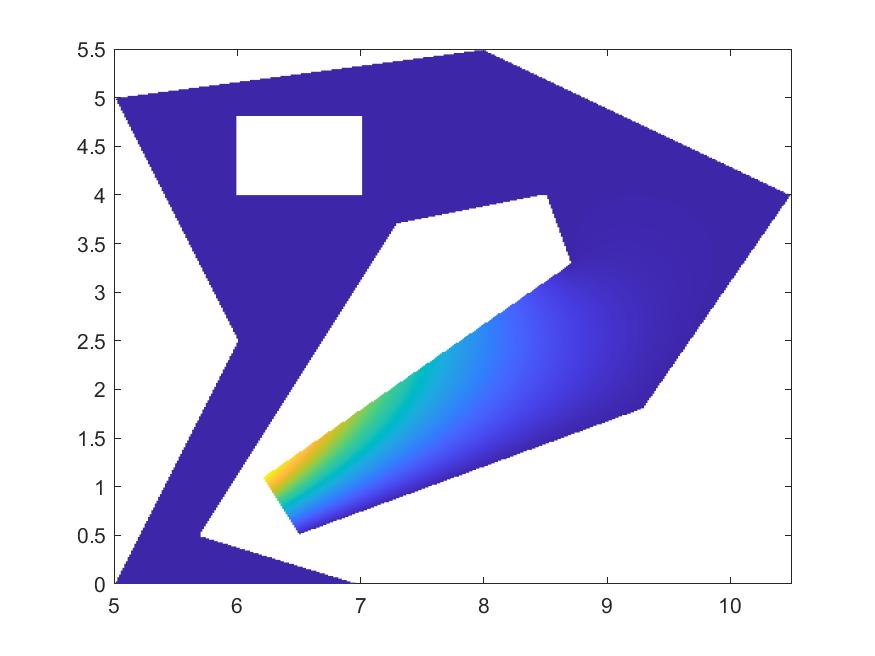}
\caption{Contour plots of GBC functions over a  domain with a hole \label{gbcholes}}
\end{figure}

\begin{figure}[h]
\centering
\includegraphics[width= 0.2\textwidth]{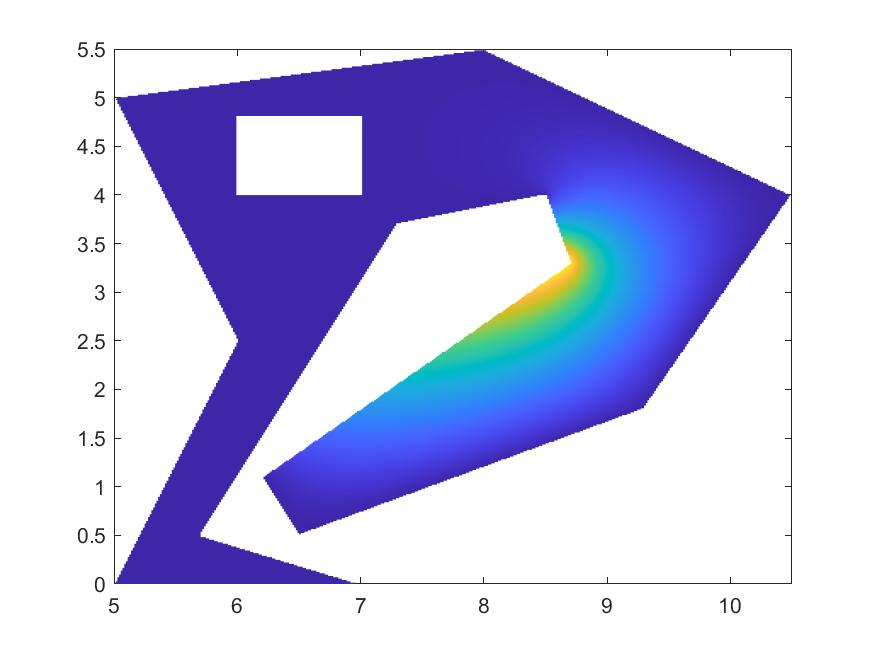}
\includegraphics[width= 0.2\textwidth]{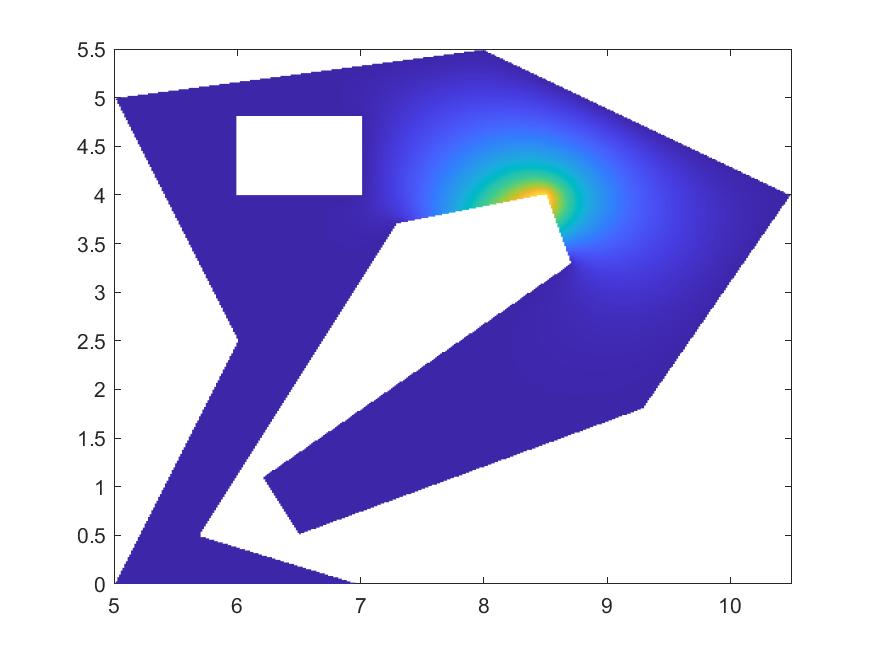}
\includegraphics[width= 0.2\textwidth]{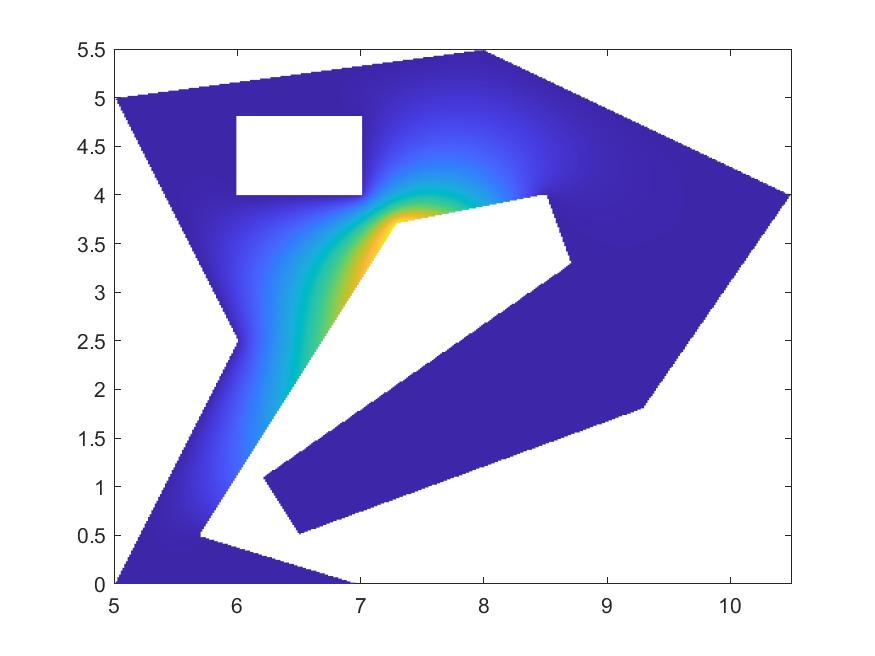}
\includegraphics[width= 0.2\textwidth]{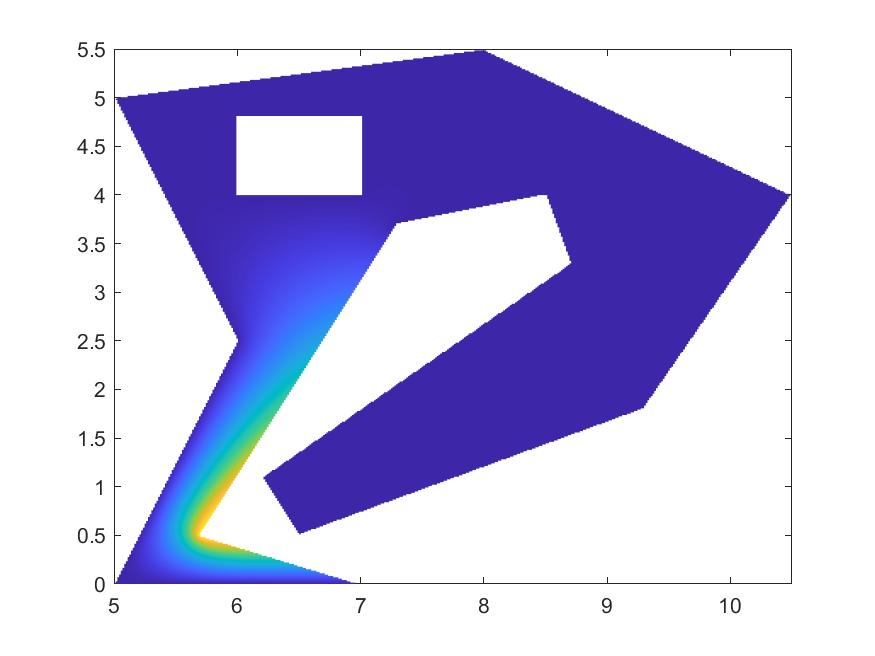}
\includegraphics[width= 0.2\textwidth]{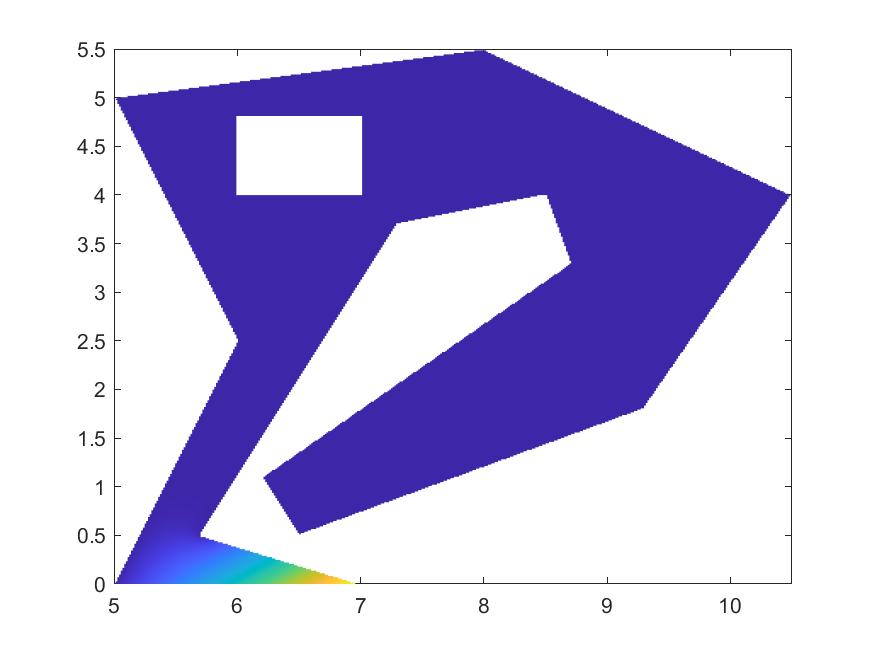}
\includegraphics[width= 0.2\textwidth]{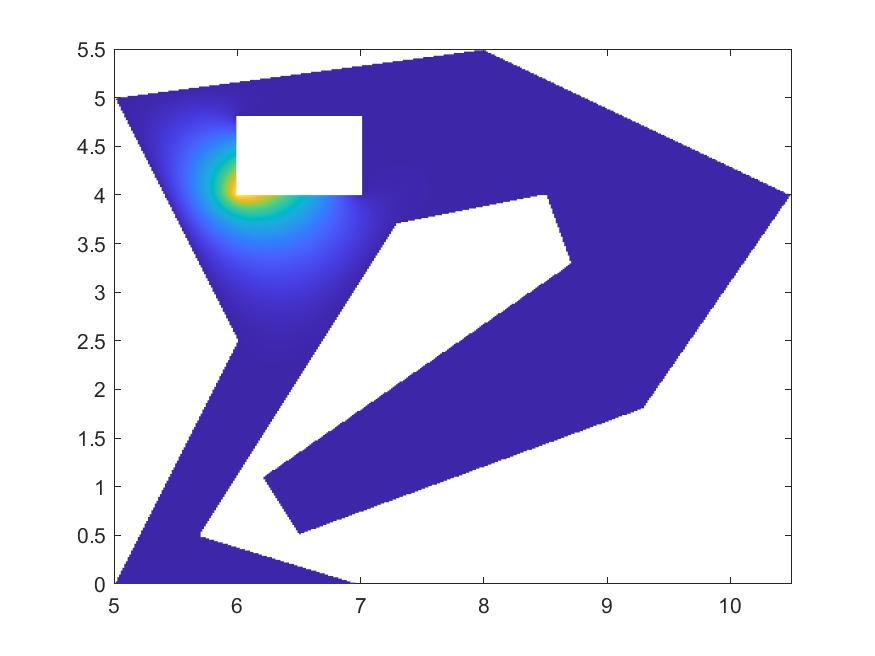}
\includegraphics[width= 0.2\textwidth]{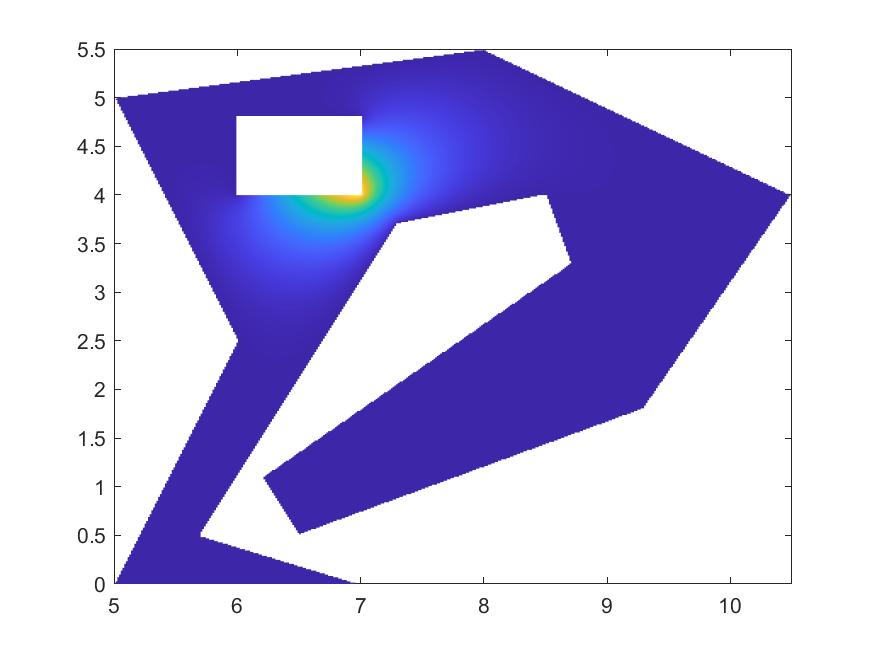}
\includegraphics[width= 0.2\textwidth]{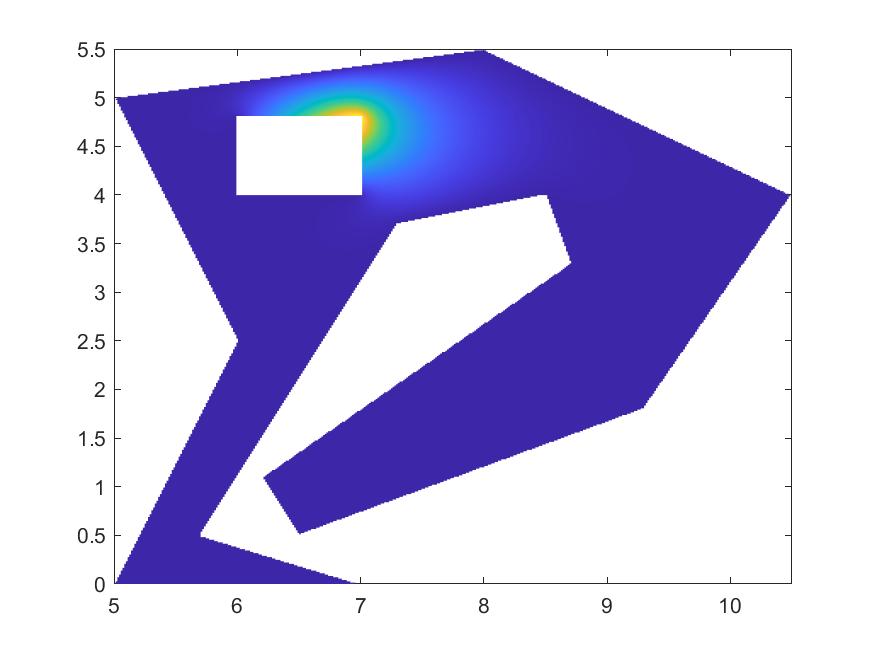}
\caption{Contour plots of more GBC functions over a  domain with a hole \label{gbcholes2}}
\end{figure}

\item 5. Although we can use 3D GBC functions to construct a deforming mapping from a 3D 
polyhedron to a convex polyhedron, say solid ball as  in \cite{SH15}, 
 it is easy to see that any extension of the proof of Theorem~\ref{harbi} to the 3D setting 
may not be easy.  This is because that the proof uses the zero property of 
holomorphic functions which has  no simple extension in the 3D setting.    
\end{itemize}

\end{document}